\documentclass[reqno]{amsart}
\usepackage{amsmath,amsfonts,amssymb,amsthm,amscd}
\newtheorem{thm}{Theorem}
\newtheorem{prop}[thm]{Proposition}
\newtheorem{cor}[thm]{Corollary}
\newtheorem{lem}[thm]{Lemma}
\newtheorem{defn}[thm]{Definition}
\def\t{{\tilde t}}

\begin{document}

\title[Conformally flat pencils and Frobenius structures]
{Conformally flat pencils of metrics, Frobenius structures and a modified Saito construction}

\author{Liana David and Ian A. B. Strachan}
\date{28$^{\rm th}$ February 2005}


\begin{abstract}
The structure of a Frobenius manifold encodes the geometry associated with a flat pencil of
metrics. However, as shown in the authors' earlier work \cite{imrn}, much of the structure
comes from the compatibility property of the pencil rather than from the flatness of the pencil itself.
In this paper conformally flat pencils of metrics are studied and examples, based on a modification
of the Saito construction, are developed.
\end{abstract}

\maketitle

\tableofcontents

\section{Introduction}

The Saito construction \cite{saito} of a flat structure on the
orbit space $\mathbb{C}^n/W\,,$ where $W$ is a Coxeter group, has
played a foundational role in many areas of mathematics. It is a
central construction in singularity theory and contains the kernel
of the definition of a Frobenius manifold, this having been done
many years before the introduction of a Frobenius manifold by
Dubrovin \cite{D}.

The initial motivation for this paper was the observation that one
may repeat the Saito construction starting with a metric of
constant non-zero sectional curvature $s$.
One easily obtains a pencil of metrics $(h,\tilde{h})$
on the orbit space
$\mathbb{C}^{n}/W$ (if $s>0$) or $\mathbb{H}^{n}\otimes\mathbb{C}/W$
(if $s<0$) of a Coexter group
$W$. The pencil $(h,\tilde{h})$ has
interesting geometric properties:
it is conformally related to the flat pencil provided by the classical
Saito construction, the metric $h$ has constant sectional curvature
$s$ and, as it turns out, the metric $\tilde{h}$ is flat.
This modified Saito construction is
developed in Section \ref{sai}. The construction may be applied,
locally, to any Frobenius manifold and this is also illustrated in
Section \ref{sai}.

Another consequence of Saito's work is that it provides a construction
of so called flat pencils of metrics. This then leads, via the results of Dubrovin and
Novikov \cite{dn}
and Magri \cite{magri}, to bi-Hamiltonian structures and
the theory of integrable systems.
The flatness of such pencils is required for the locality of the
bi-Hamiltonian structures; however one may introduce curvature - resulting in non-local
Hamiltonian operators - in such a way as to preserve the bi-Hamiltonian
property. Geometrically one requires a compatible pencil of metrics rather
than a flat pencil.

In the authors' earlier work \cite{imrn} the geometry of compatible metrics
was studied in detail - this generalizing the results of Dubrovin
\cite{dubrovin} from flat pencils of metrics
to compatible (and curved) pencils of metrics.
In Section \ref{w} we continue this study.
One way to construct
examples is to scale a
known flat pencil of metrics by a conformal factor.
This introduces
curvature but the new metrics remain compatible. The geometry of such
conformally scaled compatible pencils is studied in Section \ref{geom}
and provides a general scheme into which the modified Saito
construction of Section \ref{sai} falls.

The rest of this Section outlines some standard notations and
earlier results.

\subsection{Compatible metrics on manifolds}\label{notatii}

Let $M$ be a smooth manifold. We shall use the following notations:
${\mathcal X}(M)$ for the space of smooth vector fields on $M$;
${\mathcal E}^{1}(M)$ for the space of smooth $1$-forms
on $M$. For a pseudo-Riemannian metric $g$ on $M$,
$\nabla^{g}$ will denote
its Levi-Civita connection and $R^{g}$ its curvature.
The metric $g$ induces a metric $g^{*}$ on the cotangent bundle $T^{*}M$
of $M$. The $1$-form corresponding to
$X\in{\mathcal X}(M)$ via the pseudo-Riemannian
duality defined by $g$ will be denoted $g(X).$
Conversely,
if $\alpha\in{\mathcal E}^{1}(M)$, the corresponding vector field
will be denoted $g^{*}\alpha .$

Following \cite{imrn} we recall the basic theory of compatible metrics
on manifolds; the flat case has been treated in
\cite{dubrovin}.
Let $(g,\tilde{g})$ be an arbitrary pair of
metrics on $M$. Recall that the pair $(g,\tilde{g})$ defines
a multiplication \cite{dubrovin}
\begin{equation}\label{pr}
\alpha\circ\beta :=\nabla^{g}_{g^{*}\alpha}\beta -
\nabla^{\tilde{g}}_{g^{*}\alpha}\beta ,\quad\forall \alpha ,\beta
\in{\mathcal E}^{1}(M)
\end{equation}
on $T^{*}M$ (or on $TM$, by identifying $TM$ with $T^{*}M$ using
the metric $\tilde{g}$).
For every constant $\lambda$ we define the inverse metric
$g_{\lambda}^{*}:=g^{*}+\lambda\tilde{g}^{*}$, which, we will assume,
will always be non-degenerate and whose Levi-Civita connection
and curvature tensor will be denoted $\nabla^{\lambda}$ and $R^{\lambda}$ respectively.
The metrics $g$ and $\tilde{g}$ are almost
compatible \cite{mok} if, by definition, the relation
\begin{equation}\label{defalmost}
g_{\lambda}^{*}\nabla^{\lambda}_{X}\alpha = g^{*}\nabla^{g}_{X}\alpha
+\lambda\tilde{g}^{*}{\nabla}^{\tilde{g}}_{X}\alpha
\end{equation}
holds, for every $X\in{\mathcal X}(M)$, $\alpha\in
{\mathcal E}^{1}(M)$ and constant $\lambda .$
The almost compatibility
condition is equivalent with the vanishing of the integrability tensor
$N_{K}$ of $K:={g}^{*}\tilde{g}\in\mathrm{End}(TM)$, defined by the formula:
$$
N_{K}(X,Y)=-[KX,KY]+K[KX,Y]+K[X,KY]-K^{2}[X,Y],\quad\forall X,Y\in
{\mathcal X}(M)
$$
and implies the following two
relations:
\begin{equation}\label{aditional}
g^{*}({\nabla}^{\tilde{g}}_{\tilde{g}^{*}\gamma}\alpha
-\nabla^{g}_{\tilde{g}^{*}\gamma}\alpha ) =
\tilde{g}^{*}({\nabla}^{\tilde{g}}_{{g}^{*}\gamma}\alpha
-\nabla^{g}_{g^{*}\gamma}\alpha ) ,\quad\forall \alpha ,
\gamma\in{\mathcal E}^{1}(M)
\end{equation}
and
\begin{equation}\label{aditional1}
\tilde{g}^{*}(\alpha\circ\beta ,\gamma )=\tilde{g}^{*}(\alpha ,
\gamma\circ\beta ), \quad\forall \alpha ,\beta ,
\gamma\in{\mathcal E}^{1}(M).
\end{equation}

Recall now that two almost compatible metrics $(g,\tilde{g})$
are compatible \cite{mok} if, by definition, the relation
$$
g_{\lambda}^{*}(R^{\lambda}_{X,Y}\alpha )= g^{*}(R^{g}_{X,Y}\alpha
)+\lambda \tilde{g}^{*}({R}^{\tilde{g}}_{X,Y}\alpha )
$$
holds, for every $\alpha\in{\mathcal E}^{1}(M)$,
$X,Y\in{\mathcal X}(M)$ and constant $\lambda .$
The compatibility condition has several alternative formulations:
if the metrics $(g,\tilde{g})$ are almost compatible, then they are compatible if and
only if the relation
\begin{equation}\label{scapare}
g^{*}(\nabla^{\tilde{g}}_{X}\alpha -\nabla^{g}_{X}\alpha ,
 \nabla^{\tilde{g}}_{Y}\beta -\nabla^{g}_{Y}\beta )=
g^{*}(\nabla^{\tilde{g}}_{Y}\alpha -\nabla^{g}_{Y}\alpha ,
 \nabla^{\tilde{g}}_{X}\beta -\nabla^{g}_{X}\beta )
\end{equation}
holds for every $X,Y\in{\mathcal X}(M)$ and $\alpha ,
\beta\in{\mathcal E}^{1}(M)$,
or,
in terms of the multiplication \lq\lq$\circ$" associated to the pair $(g,\tilde{g})$,
\begin{equation}\label{sym}
(\alpha\circ\beta )\circ\gamma =(\alpha\circ \gamma )\circ\beta ,
\quad\forall\alpha ,\beta ,\gamma\in{\mathcal E}^{1}(M).
\end{equation}
If the metrics $(g,\tilde{g})$ are compatible and $R^\lambda=0$ for all
$\lambda$ then $(g,\tilde{g})$ are said to form a flat pencil of metrics
\cite{dubrovin}.\\

\subsection{The Dubrovin correspondence}

We end the Introduction by recalling the
Dubrovin correspondence \cite{dubrovin} between
flat pencils of metrics and Frobenius manifolds and its
generalizations \cite{imrn}.
We first recall the definition of a Frobenius manifold.

\begin{defn}\cite{dubrovin}
$M$ is a Frobenius manifold if a structure of a Frobenius algebra (i.e. a commutative,
associative algebra with multiplication denoted by \lq\lq$\bullet$", an
identity element \lq\lq$e$" and an
inner product \lq\lq$<,>$" satisfying the invariance condition $<a\bullet b,c>
=<a,b\bullet c>$)
is specified on the tangent space $T_pM$ at any point $p\in M$
smoothly depending on the point $p$, such that:
\begin{itemize}
\item[(i)] The invariant metric ${\tilde{g}}=<,>$ is a flat metric on
$M\,;$
\item[(ii)] The identity vector field $e$ is covariantly constant with
respect to the Levi-Civita connection
${\nabla}^{\tilde{g}}$ of the metric ${\tilde{g}}$:
\[
{\nabla}^{\tilde{g}}e=0\,;
\]
\item[(iii)] The $(4,0)$-tensor $\nabla^{\tilde{g}}(\bullet )$
defined by the formula
$$
\nabla^{\tilde{g}}(\bullet )(X,Y,Z,V):=
\tilde{g}\left( \nabla^{\tilde{g}}_{X}(\bullet )(Y,Z),V\right),\quad\forall
X,Y,Z\in {\mathcal X}(M)
$$
is symmetric in all arguments.

\item[(iv)] A vector field $E$ - the Euler vector field - must be determined on $M$ such that
\[
{\nabla}^{\tilde{g}}({\nabla}^{\tilde{g}}E)=0
\]
and that the corresponding one-parameter group of diffeomorphisms acts
by conformal transformations of the metric
and by rescalings of the multiplication \lq\lq$\bullet$".
\end{itemize}
\end{defn}

\noindent Using the flat coordinates $\{t^i\}$ of the metric $<,>$ one may
express the multiplication in terms of the derivatives of a scalar prepotential $F\,,$
\[
< \frac{\partial~}{\partial t^i} ,
\frac{\partial~}{\partial t^j} \bullet \frac{\partial~}{\partial t^k}> =
\frac{\partial^3 F}{\partial t^i \partial t^j \partial t^k}
\]
where the $t^1$-dependence of $F$ is fixed by the condition
\[
< \frac{\partial~}{\partial t^i},\frac{\partial~}{\partial t^j}>=
\frac{\partial^3 F}{\partial t^1 \partial t^i \partial t^j}\,.
\]
The associativity condition then becomes an overdetermined
partial differential equation for the prepotential $F$ known as the Witten-Dijkgraaf-Verlinde-Verlinde
equation.

Recall now that if
$(M,\bullet ,\tilde{g},E)$ is a Frobenius manifold then we can define an
inverse metric $g^{*}$ by the relation $g^{*}\tilde{g}=E\bullet .$
The metrics $(g,\tilde{g})$ form a flat pencil on the open
subset of $M$ where $E\bullet$ is an automorphism, satisfying some additional
conditions (the quasi-homogeneity conditions).
Conversely, a (regular) quasi-homogeneous flat pencil
of metrics on a manifold determine a Frobenius structure on that
manifold. This construction is known in the literature as the
Dubrovin correspondence \cite{dubrovin}.

It turns out that the key role in the Dubrovin correspondence is played
not by the flatness property of the metrics but rather by
their compatibility.
Weaker versions of the Dubrovin correspondence have been developed in
\cite{imrn}. Following \cite{imrn} we recall now a weak version of the
Dubrovin correspondence.
In general, a pair of metrics $(g,\tilde{g})$ together with a
vector field $E$ on a manifold $M$ such that the endomorphism
$T(u):=g(E)\circ u$ of $T^{*}M$ is an automorphism
(the regularity condition) determines a multiplication
$u\bullet v:=u\circ T^{-1}(v)$ on $T^{*}M$ (or on $TM$, by identifying
$TM$ with $T^{*}M$ using the metric $\tilde{g}$).
If the metrics $(g,\tilde{g})$ are compatible,
then the multiplication \lq\lq$\bullet$"
is associative, commutative,
with identity $g(E)$ on $T^{*}M$, the metrics $g$, $\tilde{g}$ are
\lq\lq$\bullet$"-invariant and
$g^{*}\tilde{g}=E\bullet .$
Moreover, if $E$ satisfies the relations
\begin{equation}\label{r}
L_{E}(\tilde{g})=D\tilde{g},\quad\nabla^{g}_{X}(E)=\frac{1-d}{2}X,\quad\forall
X\in{\mathcal X}(M),
\end{equation}
for some constants $D$ and $d$, then
$(M,\bullet ,\tilde{g},E)$ is a weak $\mathcal F$-manifold, i.e.
the following conditions are satisfied:
\begin{enumerate}
\item The metric $\tilde{g}$ and the multiplication "$\bullet$"
define a Frobenius algebra at every tangent
space of $M$.

\item The vector field $E$ - the Euler vector field - rescales
the metric $\tilde{g}$ and
the multiplication \lq\lq$\bullet$" by constants
and has an inverse $E^{-1}$ with respect to the multiplication
\lq\lq$\bullet$", which is a smooth vector field on $M$.

\item The $(4,0)$-tensor field
$\nabla^{\tilde{g}}(\bullet )$ of $M$ satisfies the symmetries:
\begin{equation}\label{simetrii}
\nabla^{\tilde{g}}(\bullet )(E,Y,Z,V)=
\nabla^{\tilde{g}}(\bullet )(Y,E,Z,V),\quad\forall Y,Z,V\in{\mathcal X}(M).
\end{equation}
\end{enumerate}

Conversely, a weak $\mathcal F$-manifold $(M,\bullet ,\tilde{g},E)$
determines a pair $(g,\tilde{g})$ of compatible metrics,
with $g$ defined by the formula $g^{*}\tilde{g}=E\bullet$, and the
Euler vector field $E$ satisfies relations $(\ref{r}).$
Therefore, there is a one to one correspondence between (regular)
compatible pencils of metrics $(g,\tilde{g})$ with a vector field $E$
satisfying relations (\ref{r}) and weak $\mathcal F$-manifolds.

Under a certain curvature
condition on the metrics $(g,\tilde{g})$ - see Theorem 23
of \cite{imrn} - the tensor $\nabla^{\tilde{g}}(\bullet )$ is symmetric in
all arguments and then
$(M,\bullet ,\tilde{g},E)$ is called an
$\mathcal F$-manifold.
Note that in this case $(M,\bullet)$ is an $F$-manifold,
i.e. the relation
\begin{equation}\label{he}
L_{X\bullet Y}(\bullet )=X\bullet L_{Y}(\bullet )+Y\bullet L_{X}(\bullet ),\quad\forall X,Y
\in{\mathcal X}(M)
\end{equation}
holds.
In fact, Hertling noticed - see Theorem 2.15 of \cite{hert} - that
if $(M,\bullet ,\tilde{g})$
satisfies the first of the three conditions mentioned above
and \lq\lq$e$" is the identity vector field of the multiplication
\lq\lq$\bullet$", then
relation (\ref{he}) together with the closeness of the coidentity
$\tilde{g}(e)$ is equivalent to
the total symmetry of the tensor $\nabla^{\tilde{g}}(\bullet ).$\\

\textbf{Remark:} The definition of (weak) $\mathcal F$-manifolds
and all the properties proved about these manifolds in \cite{imrn}
assumed that the Euler vector field rescaled the metric and the
multiplication by constants. From now on, when we refer to (weak)
$\mathcal F$-manifolds we allow the Euler vector field to rescale
the metric and the multiplication by not necessariy constant
functions. In Section \ref{w} we extend the results of \cite{imrn}
to this more general class of weak $\mathcal F$-manifolds and in
Section \ref{geom} we apply our theory to the metrics obtained by
non-constant conformal rescalings of the flat metrics of a Frobenius manifold.

\section{A modified Saito construction}\label{sai}

The motivation for considering such non-constant conformal rescalings
comes from
the following Theorem.

\begin{thm}
Let $(g,\tilde{g})$ be the flat pencil of the Saito construction
on the space of orbits $\mathbb{C}^{n}/W$ of a Coexter group $W$.
There is a metric $\tilde{h}$ with the following properties:

\begin{enumerate}

\item The metric $\tilde{h}$ is flat.

\item The metric $\tilde{h}$ is conformally related to the metric
$\tilde{g}$: $\tilde{h}=\Omega^{2}\tilde{g}$, for a smooth non-vanishing
function $\Omega .$

\item  The metric $h:=\Omega^{2}g$ has constant non-zero sectional
curvature $s$. If $s>0$ then $\tilde{h}$ is defined on $\mathbb{C}^{n}/W$. If $s<0$
then $\tilde{h}$ is defined on $\mathbb{H}^{n}\otimes\mathbb{C}/W.$
\end{enumerate}

\end{thm}

\begin{proof}
We begin with a review of the salient features of the Saito construction.
Details can be found in \cite{D}.
Recall that a Coxeter group of a real $n$-dimensional vector space
$V=\mathbb{R}^{n}$ is a finite group of linear tranformations
of $V$ generated by reflections.
Let $\{t^i\}$ be a basis of $W$-invariant
polynomials on $V$ with degrees ${\rm deg} (t^i)=d_i$, ordered so that
\[
h=d_1 > d_2 \geq \ldots \geq d_{n-1}>d_n=2\,,
\]
where $h$ is the Coxeter number of the group. The action of $W$ extends to the
complexified space $V\otimes\mathbb{C}=\mathbb{C}^{n}.$
In the Saito construction of interest is the orbit space
\[
M=\mathbb{C}^{n}/W\,.
\]
Starting with a $W$-invariant metric
\[
{\mathfrak{g}}:=\sum_{i=1}^n (dx^i)^2
\]
on $V$ one obtains a flat metric $g$ on the orbit
space $M\backslash{Discr(W)}\,,$ where
${Discr(W)}$ is the discriminant locus of irregular orbits.
What Saito showed was that there is another metric
\[
{\tilde{g}}^* := {Lie}_e (g^*)
\]
defined on the whole of $M$ which is also flat. Here $e$ is the
vector field which, in terms of the basis $\{ t^{i}\}$
of invariant polynomials,
is $\frac{\partial~}{\partial t^1}\,.$ The basis $\{ t^{i}\}$
of invariant polynomials can be chosen such that the metric
$\tilde{g}$ is anti-diagonal with constant entries:
\[
\tilde{g}_{ij}=\delta_{i+j,n+1}
\]
and is referred in this case as a Saito's basis of invariant polynomials.
Unlike $g$ which is defined only on $M\backslash{Discr(W)}\,,$ the
metric $\tilde{g}$ is defined on the
whole of $M\,.$ An important fact of the Saito construction is that
the two metrics $(g,\tilde{g})$ are
the regular flat pencil of a Frobenius structure on $M$, with Euler vector field
\begin{equation}\label{euler}
E=d_{1}t^{1}\frac{\partial}{\partial t^{1}}+\cdots +
d_{n}t^{n}\frac{\partial}{\partial t^{n}}.
\end{equation}
Suppose now that one repeats the Saito construction starting with a
metric of constant sectional
curvature, i.e. let
\begin{eqnarray*}
{\mathfrak{h}} & := & \frac{1}{\left\{c \sum_{i=1}^n (x^i)^2+d\right\}^2} {\mathfrak{g}}\,,\\
& = &\frac{1}{\left\{c\sum_{i=1}^n (x^i)^2+d\right\}^2}\sum_{i=1}^n
(dx^i)^2\,.
\end{eqnarray*}
This has constant sectional curvature $4(cd).$

Since one can take, without loss of generality, the invariant $t^n$ to be
\[
t^n= \sum_{i=1}^n (x^i)^2,
\]
the conformal factor is a function of $t^n$ alone.
Hence one obtains a new metric
\begin{eqnarray*}
{\tilde{h}}^* & := & {\rm Lie}_e (h^*)\,,\\
& = & (c t^n+d)^2 {\rm Lie}_e (g^*)\,,\\
& = & (c t^n+d)^2 {\tilde g}^*\
\end{eqnarray*}
defined on ${\mathbb{H}^n}\otimes{\mathbb{C}}/W$ for $(cd)<0$ and
on ${\mathbb{C}^{n}}/W$ for
$(cd)>0.$ In terms of the flat coordinates $\{t^i\}$ for the metric
${\tilde g}$,
\[
{\tilde{h}}_{ij} = \frac{1}{(c t^n+d)^2}\delta_{i+j,n+1}\,.
\]
Notice that the metric
$$
h:=\frac{1}{(c t^n+d )^2} g
$$
has constant sectional curvature $4(cd).$

It remains to show that
$\tilde{h}$ is flat. This may be proved using the
standard formulae for transformation of the curvature tensor under a
conformal change.
Moroever, the flat coordinates $\{{\tilde t}^i\}$ for $\tilde{h}$ can be
written down explicitly
\begin{eqnarray*}
{\tilde t}^1 & = & t^1-\frac{c}{ct^n+d} \sum_{i=2}^{n-1} t^i t^{n+1-i}\,,\\
{\tilde t}^i & = & \frac{t^i}{c t^n+d} \,,\qquad\qquad i=2\,,\ldots\,,n-1\,,\\
{\tilde t}^n & = & \frac{at^n+b}{c t^n+d}\,,\qquad\qquad ad-bc=1\,,
\end{eqnarray*}
(note that this $SL(2,\mathbb{C})$-transformation appears also in \cite{D} in a slightly
different context) giving
\[
{\tilde h}=\sum_{i=1}^n d{\tilde t}^i\, d{\tilde t}^{n+1-i}\,.
\]
The conclusion follows.
\end{proof}

The construction turns out to be quite general:

\begin{prop}\label{abra}
Suppose one has a Frobenius manifold with metrics $\tilde g$ and $g$ (or $\eta$ and $g$ respectively in
Dubrovin's notation). Consider the conformally scaled metrics
\begin{eqnarray*}
{\tilde h} & = & \Omega^2(t_1) \, {\tilde g}\,\\
h & = & \Omega^2(t_1) \, g.\,
\end{eqnarray*}
(Here $(t^{1},\cdots ,t^{n})$ are $\tilde{g}$-flat coordinates with
$\frac{\partial}{\partial t^{1}}$ being the identity vector field $e$ and
$t_1$ is dual to the identity coordinate $t^1\,,$ i.e.
$t_1=\tilde{g}_{1i}t^i$)\,.
Suppose that $\tilde h$ is flat. Then $h$ has constant sectional curvature.
\end{prop}
\begin{proof} The curvature conditions on ${\tilde h}$ translate to a
simple differential equation for the conformal factor.
Solving this gives $\Omega^{-1} = c t_1+d $ for constants $c$ and $d\,.$ This then fixes the metric
$h\,.$ Calculating its curvature (using again the standard formulae for change in the curvature tensor
under a conformal change and various properties of the Christoffel symbols of $g$ in [D]) yields the result.
\end{proof}

\noindent Note that this conformal factor satisfies the condition
$d\Omega\wedge g(E)=0$, where $E$ is the Euler vector field of the
Frobenius manifold. This follows from the following easy computation:
$g(E)=\tilde{g}(e)=\tilde{g}_{1i}dt^{i}=dt_{1}$, the functions
$\tilde{g}_{1i}$ being constant. It turns out, as Section
\ref{geom} will show, that conformally scaled metrics with this
condition have particularly attractive properties.

\section{$\mathcal F$-manifolds and compatible pencils
of metrics}\label{w}

In this Section we study the geometry of a pair of compatible
metrics together
with a vector field satisfying conditions (\ref{r}), when $D$ and $d$
are not necessarily constant.

\begin{prop}\label{cheie}
Let $(h,\tilde{h})$ be a regular pair of compatible metrics together
with a vector field $E$ on a connected manifold $M$ of dimension at least three.
Suppose that $L_{E}(h)=Dh$, $L_{E}(\tilde{h})=\tilde{D}\tilde{h}$, for
$D,\tilde{D}\in C^{\infty}(M).$
Then $E$ rescales the multiplication \lq\lq$\bullet$" associated to the pair
$(h,\tilde{h})$ and vector field $E$ if and only if
$\tilde{D}-{D}$ is constant. In this case
$L_{E}(\bullet )=(\tilde{D}-D)\bullet$ on $TM$.
\end{prop}

\begin{proof}
Recall that if $g$ is an arbitrary pseudo-Riemannian metric on a manifold
$M$ and $Z$ is a conformal vector field
with $L_{Z}(g)=pg$ for a function $p\in C^{\infty}(M)$, then
$$
L_{Z}(\nabla^{g} )_{X}(\alpha )=\frac{1}{2}[-dp(X)\alpha -
\alpha (X)dp+g^{*}(\alpha ,dp)g(X)],
$$
for every $\alpha\in{\mathcal E}^{1}(M)$ and $X\in{\mathcal X}(M).$
Applying this formula for the metrics $h$ and $\tilde{h}$ we easily get
\begin{align*}
L_{E}(\circ )(u,v)&=\frac{1}{2}[h^{*}(u,d\lambda )v+h^{*}(u,v)d\lambda
+h^{*}(v,dD)u-\tilde{h}^{*}(v,d\tilde{D})\tilde{h}h^{*}(u)]\\
&-Du\circ v,
\end{align*}
where \lq\lq$\circ$"
is the multiplication (\ref{pr}) determined by the pair of metrics
$(h,\tilde{h})$, $u,v\in{\mathcal E}^{1}(M)$ and $\lambda :=\tilde{D}-D$.
Let $T$ be the automorphism of $T^{*}M$ defined by the formula
$T(u)=h(E)\circ u$. From the above relation we easily see that
$$
L_{E}(T)(u)=\frac{1}{2}[E(\lambda )u+u(E)d\lambda +
h^{*}(u,dD )h(E)-\tilde{h}^{*}(u,d\tilde{D})\tilde{h}(E)].
$$

Recall now that the multiplications \lq\lq$\circ$" and \lq\lq$\bullet$"
on $T^{*}M$ are related by the formula $u\bullet T(v)=u\circ v.$ Taking the derivative
with respect to $E$ of this formula we easily see that
\begin{align*}
L_{E}(\bullet )(u,Tv)&+\frac{1}{2}u\bullet [E(\lambda )v+v(E)d\lambda +h^{*}(v,dD)h(E)-
\tilde{h}^{*}(v,d\tilde{D})\tilde{h}(E)]\\
&=\frac{1}{2}[h^{*}(u,d\lambda )v+h^{*}(u,v)d\lambda +h^{*}(v,dD)u
-\tilde{h}^{*}(v,d\tilde{D})\tilde{h}h^{*}(u)]\\
&-Du\bullet T(v),
\end{align*}
for every $u,v\in{\mathcal E}^{1}(M).$
Since $h^{*}\tilde{h}=E\bullet $
(the metrics $(h,\tilde{h})$ being compatible), an easy argument shows
that
$u\bullet\tilde{h}(E)=\tilde{h}h^{*}(u)$.
The compatibility of $(h,\tilde{h})$
also implies, as mentioned in Section \ref{notatii}, that $h(E)$ is the identity
of the multiplication "$\bullet$" on $T^{*}M.$
It follows that
\begin{align*}
L_{E}(\bullet )(u,T(v))&=\frac{1}{2}[h^{*}(u,d\lambda )v+h^{*}(u,v)d\lambda ]\\
&-\frac{1}{2}u\bullet [E(\lambda )v+v(E)d\lambda ]\\
&-Du\bullet T(v),
\end{align*}
for every $u,v\in{\mathcal E}^{1}(M).$
{}From this relation it is easy to see that $E$ rescales the multiplication
\lq\lq$\bullet$"
if and only if
for every $u,v\in{\mathcal E}^{1}(M)$, the equality
\begin{equation}\label{lung}
u\bullet [E(\lambda )v+v(E)d\lambda ]=h^{*}(u,d\lambda )v+h^{*}(u,v)d\lambda
\end{equation}
holds and in this case $L_{E}(\bullet )=-D\bullet$ on $T^{*}M$, or
$L_{E}(\bullet )=\lambda\bullet$ on $TM$. We will show that relation
(\ref{lung}) holds only if $\lambda$ is constant. Indeed, if
in relation (\ref{lung}) we take $u$ and $v$ annihilating $E$, then
we get
$$
h^{*}(u,d\lambda )v=h^{*}(v,d\lambda )u
$$
which can hold only if
$d\lambda =\mu h(E)$ for a function $\mu\in C^{\infty}(M)$,
since the dimension
of $M$ is at least three (and hence the annihilator of $E$ in
$T^{*}M$ is of dimension at least two). Relation (\ref{lung}) then becomes
$$
\mu [h(E,E)v+v(E)h(E)]\bullet u=\mu [u(E)v+h^{*}(u,v)h(E)]
$$
which in turn implies that
$\mu [v(E)u-u(E)v]$ is symmetric in $u$ and $v$, for every $u,v\in
{\mathcal E}^{1}(M).$ This can happen
only when $\mu$ is identically zero or $\lambda =\tilde{D}-D$ is constant
($M$ being connected).
\end{proof}

\textbf{Remark:} Note that relation (\ref{lung}) does not
imply that $\lambda$ is constant in dimension two. An easy argument shows
that relation (\ref{lung})
imposes that in two dimensions the
multiplication \lq\lq$\bullet$"
is of the form
$$
d\lambda\bullet d\lambda =0;\quad d\lambda\bullet h(E)=d\lambda ;\quad
h(E)\bullet h(E)=h(E)
$$
when $\lambda$ is non-constant. Proposition \ref{cheie} does not hold in
dimension two: consider for example the inverse metrics
\begin{align*}
\tilde{h}^{*}&=
f\left(\frac{\partial}{\partial x}
\otimes\frac{\partial}{\partial y}
+\frac{\partial}{\partial y}
\otimes\frac{\partial}{\partial x}
+\frac{\partial}{\partial y}\otimes\frac{\partial}{\partial y}\right)\\
{h}^{*}&=x\left(\frac{\partial}{\partial x}
\otimes\frac{\partial}{\partial y}
+\frac{\partial}{\partial y}
\otimes\frac{\partial}{\partial x}\right)
+y\frac{\partial}{\partial y}\otimes\frac{\partial}{\partial y}
\end{align*}
together with the vector field
$$
E=x\frac{\partial}{\partial x}+y\frac{\partial}{\partial y},
$$
with $f$ smooth, non-vanishing, depending only on $x$,
such that $\frac{xf^{\prime}(x)}{f(x)}$ is non-constant
on a connected open subset $M$ of
$$
\{ (x,y)\in\mathbb{R}^{2}:\quad
\frac{xf^{\prime}(x)}{f(x)}\neq\frac{1}{2},\quad x\neq 0\} .
$$
The hypothesis of Proposition \ref{cheie} is satisfied
is satisfied on $M$, with
$$
L_{E}(h)=h;\quad L_{E}(\tilde{h})
=\left( 2-\frac{xf^{\prime}(x)}{f(x)}\right)
\tilde{h}.
$$
The associated multiplication \lq\lq$\bullet$" on $T^{*}M$ has the
expression:
$$
dx\bullet dx=0;\quad dx\bullet dy=dx;\quad dy\bullet dy=dy
$$
and is preserved by $E$:
$L_{E}(\bullet) = \left( 1-\frac{xf^{\prime}(x)}{f(x)}\right)
\bullet$ on $TM.$\\

We return now to higher dimensions and we note the following consequence
of Proposition \ref{cheie}.

\begin{cor}\label{dim}
Let $(M,\bullet ,\tilde{h},E)$ be a connected weak
$\mathcal F$-manifold of dimension
at least three. Then
the Euler vector field $E$ rescales the multiplication \lq\lq$\bullet$" by
a constant.
\end{cor}

\begin{proof}
Let $h$ be the metric on $M$ defined by the relation
$h^{*}\tilde{h}=E\bullet .$
Suppose that $L_{E}(\tilde{h})=\tilde{D}\tilde{h}$ and $L_{E}(\bullet )=
k\bullet$, for $\tilde{D},k\in C^{\infty}(M).$
Then $L_{E}(h)=(\tilde{D}-k)h.$
The same argument used in the proof of Theorem 17 and Proposition 19
of \cite{imrn} shows that $(h,\tilde{h})$ are compatible and that
the multiplication \lq\lq$\bullet$" on $T^{*}M$ (identified with $TM$
using the metric $\tilde{h}$) is related to the multiplication \lq\lq$\circ$"
determined by the pair $(h,\tilde{h})$ by the relation:
$u\bullet T(v)=u\circ v$, where $T$ is the endomorphism of $T^{*}M$
defined by $T(u):=h(E)\circ u.$
Even if $T$ is not necessarily an automorphism, the argument
used in the proof of Proposition \ref{cheie} still holds and implies
our conclusion.

\end{proof}

In order to simplify terminology, we introduce the following
definition (which generalizes Definition 14 of \cite{imrn}).

\begin{defn}\label{defini}
Let $(h,\tilde{h})$ be a compatible pair of metrics and $E$ a
vector field on a manifold $M$. The pair $(h,\tilde{h})$ is a weak
quasi-homogeneous pencil with Euler vector field $E$ if
the following conditions are satisfied:

\begin{enumerate}

\item $L_{E}(\tilde{h})=\tilde{D}\tilde{h}$;
$\nabla^{h}(E)=\frac{D}{2}\mathrm{Id}$, where \lq\lq$\mathrm{Id}$" is the
identity endomorphism of $TM$ and $D,\tilde{D}\in C^{\infty}(M).$

\item The difference $\tilde{D}-{D}$ is constant.

\end{enumerate}
The weak quasi-homogeneous pair $(h,\tilde{h})$ is regular if
the endomorphism $T(u):=h(E)\circ u$ of $T^{*}M$ is an automorphism.
(Here \lq\lq$\circ$" is the multiplication (\ref{pr})
associated to the pair of metrics $(h,\tilde{h})$).

\end{defn}

The correspondence between weak quasi-homogeneous
pencils of metrics and weak $\mathcal F$-manifolds can be stated as follows:

\begin{thm}

\begin{enumerate}

\item Let $(h,\tilde{h})$ be a regular weak quasi-homogeneous pencil of metrics
with Euler vector field $E$ on a manifold $M$.
Let \lq\lq$\bullet$" be the multiplication on $TM$ associated to the
pair of metrics $(h,\tilde{h})$ and vector field $E$.
Then $(M,\bullet ,\tilde{h},E)$ is a weak $\mathcal F$-manifold.

\item Conversely, let $(M,\bullet ,\tilde{h},E)$ be a connected weak
$\mathcal F$-manifold of dimension at least three.
Define the metric $h$ on $M$ by the formula $h^{*}\tilde{h}=E\bullet .$
Then $(h,\tilde{h})$
is a weak quasi-homogeneous pencil with Euler vector field $E$.
\end{enumerate}
\end{thm}

\begin{proof}
The proof follows the same steps and is totally similar
to the proofs of Theorem 17 and Theorem 20 of \cite{imrn}.
Note that in the second statement of the Theorem we have restricted
to the case when the manifold $M$ is connected and of dimension at least three.
These additional conditions insure that the pair $(h,\tilde{h})$
satisfies the second condition of Definition \ref{defini}.

\end{proof}

The following Theorem and its Corollary
generalize the results from Section 6 of \cite{imrn}.

\begin{thm}\label{gen}
Let $(M,\bullet ,\tilde{h},E)$ be a weak $\mathcal F$-manifold with
$L_{E}(\tilde{h})=\tilde{D}\tilde{h}$,
$L_{E}(\bullet )=k\bullet$, where $k$ is constant
and $\tilde{D}\in C^{\infty}(M).$
Let $h$ be the metric on $M$ defined by the relation $h^{*}\tilde{h}=E\bullet .$
Consider the multiplication \lq\lq$\bullet$" also on $T^{*}M$, by
identifying $TM$ and $T^{*}M$ using the metric $\tilde{h}.$
Then $(M,\bullet ,\tilde{h},E)$ is
an $\mathcal F$-manifold if and only if the equality
\begin{align*}
R^{h}_{X,Y}(\alpha )=R^{\tilde{h}}_{X,Y}(\alpha )&+
\left( -R^{\tilde{h}}_{X,E}(\alpha )+\frac{1}{2}\alpha (X)d\tilde{D}\right)
\bullet\tilde{h}
(E^{-1}\bullet Y)\\
&-\left( -R^{\tilde{h}}_{Y,E}(\alpha )+\frac{1}{2}\alpha (Y)d\tilde{D}\right)
\bullet\tilde{h}
(E^{-1}\bullet X)\\
\end{align*}
holds, for every $X,Y\in{\mathcal X}(M)$ and $\alpha\in{\mathcal E}^{1}(M).$
\end{thm}

\begin{proof}
The argument is similar to the one employed in the proof of
Theorem 23 of \cite{imrn}. The only difference from the case studied in \cite{imrn} is that $\tilde{D}$
can be non-constant and then
$$
\nabla^{\tilde{h}}_{Y}(\nabla^{\tilde{h}} E)_{X}=R^{\tilde{h}}_{Y,E}(X)+
\frac{1}{2}[-\tilde{h}(X,Y)d(\tilde{D})+
Y(\tilde{D})\tilde{h}(X)+
X(\tilde{D})\tilde{h}(Y)]
$$
for every $X,Y\in{\mathcal X}(M)$, which is the analogue of
Lemma 22 of \cite{imrn} and can be proved in the same way in this more
general context.
\end{proof}

\begin{cor}\label{cu}
Consider the set-up of Proposition \ref{gen} and suppose that
$\tilde{h}$ is flat.
Then $(M,\bullet ,\tilde{h},E)$ is
an $\mathcal F$-manifold if and only if $h$ has constant
sectional curvature $s$ and $d\tilde{D}=-2sh(E)$.
\end{cor}

\begin{proof}
{}From Theorem \ref{gen} and the flatness of $\tilde{h}$
we know that $(M,\bullet ,\tilde{h}, E)$
is an $\mathcal F$-manifold if and only if the curvature
$R^{h}$ of $h$ has the following expression:
\begin{equation}\label{echiv}
R^{h}_{X,Y}(\alpha )=\frac{1}{2}d\tilde{D}\bullet\left(\alpha (X)
\tilde{h}(E^{-1}\bullet Y)
-\alpha (Y)\tilde{h}(E^{-1}\bullet X)\right) ,
\end{equation}
for every $X,Y\in{\mathcal X}(M)$ and $\alpha\in{\mathcal E}^{1}(M).$
It is clear now that if $h$ has constant sectional curvature $s$ and
$d\tilde{D}=-2sh(E)$, then relation (\ref{echiv}) is satisfied.
Conversely, suppose that $(M,\bullet ,\tilde{h},E)$ is an
$\mathcal F$-manifold, so that relation (\ref{echiv}) is satisfied.
Then
\begin{equation}\label{e}
h\left( R^{h}_{X,Y}Z,V\right) =
\frac{1}{2}\left(h(X,Z)d\tilde{D}(V\bullet Y\bullet E^{-1})-
h(Y,Z)d\tilde{D}(V\bullet X\bullet E^{-1})\right),
\end{equation}
for every $X,Y,Z,V\in{\mathcal X}(M).$
On the other hand, since
$$
h\left( R^{h}_{X,Y}Y,X\right) =-h\left( R^{h}_{X,Y}X,
Y\right),\quad\forall X,Y\in{\mathcal X}(M)
$$
we easily get
$$
h(X,X)(d\tilde{D})(Y^{2}\bullet E^{-1})=h(Y,Y)(d\tilde{D})(X^{2}\bullet E^{-1})
$$
or
$$
h(X,T)(d\tilde{D})(Y\bullet S\bullet E^{-1})=h(Y,S)(d\tilde{D})(X\bullet T\bullet E^{-1})
,\quad\forall X,Y,S,T\in{\mathcal X}(M).
$$
It follows that $h(E)\wedge d\tilde{D}=0$ (let $S=T:=E$ in the above relation)
or $d\tilde{D}=-2s h(E)$, for a function $s\in C^{\infty}(M).$
{}From relation (\ref{e}) we deduce that $s$ is constant and
$h$ has constant sectional curvature $s$.

\end{proof}

\section{The geometry of conformally scaled compatible pencils}\label{geom}

In this Section we fix a pair of metrics $(g,\tilde{g})$ on a manifold
$M$. The following Lemma will be relevant in our calculations.

\begin{lem}\label{ac}
Suppose that the metrics $(g,\tilde{g})$ are almost compatible.
Then, for every
$X,Y\in{\mathcal X}(M)$ and $\alpha\in{\mathcal E}^{1}(M)$ the relation
$$
g^{*}\left(\nabla^{\tilde{g}}_{X}\alpha -\nabla^{g}_{X}\alpha,
\tilde{g}(Y)\right)
=g^{*}\left( \nabla^{\tilde{g}}_{Y}\alpha -\nabla^{g}_{Y}\alpha ,
\tilde{g}(X)\right)
$$
holds.
\end{lem}

\begin{proof}
Let $X:=\tilde{g}^{*}(\gamma )$ and $Y:=\tilde{g}^{*}(\delta )$, for
$\gamma ,\delta\in{\mathcal E}^{1}(M).$ Then
\begin{align*}
g^{*}\left(\nabla^{\tilde{g}}_{X}\alpha
-\nabla^{g}_{X}\alpha ,\tilde{g}(Y)\right)&=
\delta\left( g^{*}(\nabla^{\tilde{g}}_{\tilde{g}^{*}\gamma}\alpha
-\nabla^{g}_{\tilde{g}^{*}\gamma}\alpha )\right)\\
&=\delta\left(\tilde{g}^{*}(\nabla^{\tilde{g}}_{g^{*}\gamma}\alpha -
\nabla^{g}_{g^{*}\gamma}\alpha )\right)\\
&=\tilde{g}^{*}(\gamma\circ\alpha ,\delta )=\tilde{g}^{*}
(\gamma ,\delta\circ\alpha)\\
&=g^{*}\left( \nabla^{\tilde{g}}_{Y}\alpha -\nabla^{g}_{Y}\alpha
,\tilde{g}(X)\right),
\end{align*}
where
\lq\lq$\circ$" is the multiplication (\ref{pr}) associated to the pair
$(h,\tilde{h})$ and
we have used
relations (\ref{aditional}) and (\ref{aditional1}).
\end{proof}

As a consequence of Lemma \ref{ac} we deduce that the compatibility
property of two metrics is conformal invariant:

\begin{prop}\label{comp} Suppose that the metrics
$(g,\tilde{g})$ are compatible and let
$\Omega\in C^{\infty}(M)$, non-vanishing.
Then the metrics
$(h:=\Omega^{2}g,\tilde{h}:=\Omega^{2}\tilde{g})$ are also compatible.
\end{prop}

\begin{proof}

It is obvious that the metrics $h$ and $\tilde{h}$ are almost compatible,
since
${h}^{*}\tilde{h}={g}^{*}\tilde{g}$ (and hence the integrability tensor of
${h}^{*}\tilde{h}$ is identically zero). In order to show the compatibility
of $(h,\tilde{h})$,
we first notice that
\begin{align*}
\nabla^{h}_{X}\alpha &=\nabla^{g}_{X}\alpha -\frac{d\Omega}{\Omega}(X)
\alpha -\alpha (X)\frac{d\Omega}{\Omega}+
g^{*}\left(\alpha ,\frac{d\Omega}{\Omega}\right)g(X)\\
\nabla^{\tilde{h}}_{X}\alpha &=\nabla^{\tilde{g}}_{X}\alpha -
\frac{d\Omega}{\Omega}(X)\alpha -
\alpha (X)\frac{d\Omega}{\Omega}+\tilde{g}^{*}\left(
\alpha ,\frac{d\Omega}{\Omega}\right)\tilde{g}(X),
\end{align*}
for every $X\in{\mathcal X}(M)$ and $\alpha\in\mathcal E^{1}(M)$, from where we deduce that
\begin{equation}\label{varc}
\nabla^{\tilde{h}}_{X}\alpha -\nabla^{h}_{X}\alpha
=\nabla^{\tilde{g}}_{X}\alpha
-\nabla^{g}_{X}\alpha +\tilde{g}^{*}\left(\alpha ,\frac{d\Omega}{\Omega}
\right)\tilde{g}(X)- g^{*}\left(\alpha ,\frac{d\Omega}{\Omega}\right)g(X).
\end{equation}
To prove the compatibility of the metrics $(h,\tilde{h})$ we shall verify
relation (\ref{scapare}).
Notice that, since $h^{*}=\Omega^{-2}g^{*}$,
we need to show that the relation
\begin{equation}\label{au}
g^{*}(\nabla^{\tilde{h}}_{X}\alpha -\nabla^{h}_{X}\alpha ,
 \nabla^{\tilde{h}}_{Y}\beta -\nabla^{h}_{Y}\beta )=
g^{*}(\nabla^{\tilde{h}}_{Y}\alpha -\nabla^{h}_{Y}\alpha ,
 \nabla^{\tilde{h}}_{X}\beta -\nabla^{h}_{X}\beta )
\end{equation}
holds, for every $X,Y\in{\mathcal X}(M)$ and
$\alpha ,\beta\in{\mathcal E}^{1}(M)$. Using the
compatibility of the metrics $(g,\tilde{g})$ and relation (\ref{varc}),
we easily see that
relation (\ref{au}) is equivalent with
\begin{align*}
\tilde{g}^{*}\left(\beta ,\frac{d\Omega}{\Omega}\right)
&[g^{*}\left(\nabla^{\tilde{g}}_{X}\alpha -\nabla^{g}_{X}
\alpha,\tilde{g}(Y)\right)-g^{*}\left(\nabla^{\tilde{g}}_{Y}\alpha -
\nabla^{g}_{Y}\alpha ,\tilde{g}(X)\right) ]+\\
\tilde{g}^{*}\left(\alpha ,\frac{d\Omega}{\Omega}\right)
&[g^{*}\left(\nabla^{\tilde{g}}_{Y}\beta
-\nabla^{g}_{Y}
\beta,\tilde{g}(X)\right)-g^{*}\left(\nabla^{\tilde{g}}_{X}\beta -
\nabla^{g}_{X}\beta ,\tilde{g}(Y)\right) ]=0,\\
\end{align*}
which is obviously true from Lemma \ref{ac}.
\end{proof}

For the rest of this Section we suppose that the metrics
$(g,\tilde{g})$ are
the regular flat metrics of a Frobenius manifold $(M,\bullet ,\tilde{g},E)$.
We study the geometry of the pair of scaled
metrics $(h:=\Omega^{2}g,
\tilde{h}:=\Omega^{2}\tilde{g})$ together with the vector field $E$.
We restrict to the case when the scaled pair is regular and we denote by
\lq\lq$\bullet_{h}$" the associated multiplication on $TM$ or $T^{*}M.$
Recall that the multiplication \lq\lq$\bullet_{g}$" on $TM$ associated to the pair
of metrics $(g,\tilde{g})$ together with $E$ coincides with the
multiplication \lq\lq$\bullet$" of the Frobenius manifold $(M,\bullet
,\tilde{g},E).$

\begin{prop}\label{mul}
The multiplications \lq\lq$\bullet_{h}$" and \lq\lq$\bullet_{g}$"
coincide on $TM.$
\end{prop}

\begin{proof}

{}From relation (\ref{varc}) we easily see that
the multiplications \lq\lq$\circ_{h}$" and \lq\lq$\circ_{g}$"
associated to the pair of metrics $(h,\tilde{h})$ and
$(g,\tilde{g})$ respectively
are related by the formula
\begin{equation}\label{aux1}
\alpha\circ_{h}\beta =\Omega^{-2}[\alpha\circ_{g}\beta +g^{*}\left(\beta ,
\frac{d\Omega}{\Omega}\right)\alpha -\tilde{g}^{*}\left(\beta ,
\frac{d\Omega}{\Omega}\right)\tilde{g}g^{*}(\alpha )],
\end{equation}
for every $\alpha ,\beta\in{\mathcal E}^{1}(M)$.
Define the automorphisms
$T(\alpha ):=g(E)\circ_{g}\alpha$ and $\tilde{T}(\alpha ):=h(E)\circ_{h}
\alpha $ of $T^{*}M.$ Relation (\ref{varc}) also implies that
\begin{equation}\label{aux2}
\tilde{T}(\alpha )=T(\alpha )+g^{*}\left( \alpha,
\frac{d\Omega}{\Omega}\right) g(E)
-\tilde{g}^{*}\left( \alpha,\frac{d\Omega}{\Omega}\right)\tilde{g}(E).
\end{equation}
Since $\alpha\circ_{h}\beta =\alpha\bullet_{h}\tilde{T}(\beta)$ and
similarly $\alpha\circ_{g}\beta =\alpha\bullet_{g}{T}(\beta)$ we
deduce from (\ref{aux1}) and (\ref{aux2})
that the relation
\begin{align*}
\alpha\bullet_{h}&[T(\beta )+g^{*}\left(\beta ,
\frac{d\Omega}{\Omega}\right) g(E)
-\tilde{g}^{*}\left( \beta ,\frac{d\Omega}{\Omega}\right)\tilde{g}(E)]\\
&=\Omega^{-2}[\alpha\bullet_{g} T(\beta )+g^{*}\left(\beta,
\frac{d\Omega}{\Omega}\right)\alpha -\tilde{g}^{*}\left(\beta ,
\frac{d\Omega}{\Omega}\right)\tilde{g}g^{*}(\alpha )]
\end{align*}
holds, for every $\alpha ,\beta\in{\mathcal E}^{1}(M)$.
As in the
proof of Proposition \ref{cheie}, $\alpha\bullet_{h}\tilde{g}(E)=\Omega^{-2}
\tilde{h}h^{*}(\alpha)$ and $\alpha\bullet_{h}g(E)=\Omega^{-2}\alpha$
(the metrics $(h,\tilde{h})$ being compatible).
It follows that
$\bullet_{h}=\Omega^{-2}\bullet_{g}$ on $T^{*}M$, or
$\bullet_{h}=\bullet_{g}$ on $TM.$

\end{proof}


\begin{prop}\label{f}
The following statements are equivalent:
\begin{enumerate}

\item the pair $(h,\tilde{h})$ is weak quasi-homogeneous
with Euler vector field $E.$

\item $g(E)\wedge d\Omega =0.$

\item $(M,\bullet ,\tilde{h},E)$ is an $\mathcal F$-manifold.

\item $(M,\bullet ,\tilde{h},E)$ is a weak $\mathcal F$-manifold.

\end{enumerate}
\end{prop}

\begin{proof}

Before proving the equivalence of the statements, we make some preliminary
remarks.
Since $(g,\tilde{g})$ are the flat metrics of a
Frobenius manifold, $L_{E}(g)=(1-d)g$ and $L_{E}(\tilde{g})=D\tilde{g}$ for
some constants $D$ and $d$. It follows that
$$
L_{E}(h)=\left(1-d+\frac{2E(\Omega)}{\Omega}\right)h,
\quad L_{E}(\tilde{h})=\left(D+\frac{2E(\Omega )}{\Omega}\right)\tilde{h}.
$$
Also,
\begin{align*}
\nabla^{h}_{X}(E)&=\nabla_{X}^{g}(E)+\frac{d\Omega}{\Omega}(X)E+
\frac{E(\Omega )}{\Omega}X-g(X,E)g^{*}\left(\frac{d\Omega}{\Omega}\right)\\
&=\left(\frac{1-d}{2}+\frac{E(\Omega )}{\Omega}\right) X-
\left( E\wedge{g}^{*}\left(\frac{d\Omega}{\Omega}\right)\right)
({g}(X)).
\end{align*}
Moreover, Proposition \ref{comp} implies that the metrics $(h,\tilde{h})$ are
compatible.
The equivalence $1\Longleftrightarrow 2$ clearly follows from these facts.
The equivalence $2\Longleftrightarrow 3$ follows from Hertling's
observation mentioned at the end of
Section \ref{notatii}: indeed, condition $(2)$ means that the coidentity
$\tilde{h}(e)=\Omega^{2}g(E)$ is closed (note that the $1$-form
$g(E)$, being equal to $\tilde{g}(e)$, is closed because
$e$ is $\nabla^{\tilde{g}}$-parallel).
To prove the equivalence $3\Longleftrightarrow 4$
we notice that, since $(M,\bullet ,\tilde{g},E)$ is an
$\mathcal F$-manifold, the $(3,1)$-tensor field
$\nabla^{\tilde{h}}(\bullet )$ satisfies
the relation
\begin{align*}
\nabla^{\tilde{h}}_{X}(\bullet )(Y,Z)-\nabla^{\tilde{h}}_{Y}(\bullet )(X,Z)
&=\left(\tilde{g}^{*}\left(\frac{d\Omega}{\Omega}\right)\wedge e\right)
\left(\tilde{g}(Y\bullet Z)\right)\bullet X\\
&-\left(\tilde{g}^{*}\left(\frac{d\Omega}{\Omega}\right)\wedge e\right)
\left(\tilde{g}(X\bullet Z)\right)\bullet Y,\\
\end{align*}
for every $X,Y,Z\in{\mathcal X}(M)$.
Suppose now that $(M,\bullet ,\tilde{h},E)$ is a weak $\mathcal F$-manifold.
The symmetry
$$
\nabla^{\tilde{h}}_{E}(\bullet )(Y,Z)=\nabla^{\tilde{h}}_{Y}(\bullet )(E,Z),
\quad\forall Y,Z\in{\mathcal X}(M)
$$
of the $(3,1)$-tensor field $\nabla^{\tilde{h}}(\bullet )$ becomes,
after replacing $Z$ with $E^{-1}\bullet Z$, the relation:
$$
\left(\tilde{g}^{*}\left(\frac{d\Omega}{\Omega}\right)\wedge e\right)
\left(\tilde{g}(Z)\right)\bullet Y=
\left(\tilde{g}^{*}\left(\frac{d\Omega}{\Omega}\right)\wedge e\right)
\left(\tilde{g}(Y\bullet Z\bullet E^{-1})\right)\bullet E.
$$
It is clear now that if $(M,\bullet ,\tilde{h},E)$ is a weak
$\mathcal F$-manifold, then
$$
\nabla^{\tilde{h}}_{X}(\bullet )(Y,Z)=\nabla^{\tilde{h}}_{Y}(\bullet )(X,Z),
\quad\forall X,Y,Z\in {\mathcal X}(M)
$$
which implies that $(M,\bullet ,\tilde{h},E)$ is an $\mathcal F$-manifold
(the symmetry of the $(4,0)$-tensor field  $\nabla^{\tilde{h}}(\bullet )$
in the last three arguments is a
consequence of the fact that $\tilde{h}$
is \lq\lq$\bullet$"-invariant and of the commutativity of \lq\lq$\bullet$").
The equivalence
$3\Longleftrightarrow 4$ follows.

\end{proof}

Note that Hertling's observation used in the proof
of Proposition \ref{f} together with Corollary \ref{cu}
provide a different view-point of Proposition \ref{abra}.

\begin{cor}
Let $(g,\tilde{g})$ be the flat metrics of a Frobenius
manifold $(M,\bullet ,\tilde{g},E).$ Let $\Omega\in C^{\infty}(M)$
non-vanishing which satisfies $d\Omega\wedge g(E)=0$. Consider the scaled
metrics $(h:=\Omega^{2}g,\tilde{h}:=\Omega^{2}\tilde{g}).$
If $\tilde{h}$ is flat, then $h$ has constant sectional curvature.
\end{cor}

\begin{proof}
The condition $d\Omega\wedge g(E)=0$ implies, using Hertling's observation,
that $(M,\bullet ,\tilde{h},E)$ is an $\mathcal F$-manifold.
The conclusion follows from Corollary \ref{cu}
(since $h^{*}\tilde{h}=E\bullet$ and $\tilde{h}$ is flat).

\end{proof}

\section{The modified Saito construction revisited}

We return now to the modified Saito construction
described in Section 2, summarizing the various results in the
following Theorem.

\begin{thm} Let $(g,\tilde{g})$ be the flat metrics of the Frobenius
structure $(M,\bullet ,\tilde{g}, E)$ on the space of orbits
$M=\mathbb{C}^{n}/W$
of a Coexter group $W$.
Let $\{ t^{i}\}$ be a Saito basis of $W$-invariant polynomials
with $\mathrm{deg}(t^{i})=d_{i}$,
in which the metric $\tilde{g}$ is anti-diagonal,
the identity vector field $e$ is $\frac{\partial}{\partial t^{1}}$
and the Euler vector field $E$ has the expression (\ref{euler}).
Consider the pair of scaled metrics $(h:=\Omega^{2}g,\tilde{h}:
=\Omega^{2}\tilde{g})$,
where $\Omega\in C^{\infty}(M_{0})$ is non-vanishing
on an open subset $M_{0}$ of $M.$
The following facts hold:

\begin{enumerate}

\item The metrics $(h,\tilde{h})$ are compatible on $M_{0}$.

\item The metrics $(h,\tilde{h})$ together with the Euler vector field
$E$ is a weak quasi-homogeneous pencil on $M_{0}$ if and only
if $\Omega$ depends
only on the last coordinate $t^{n}.$
If $\Omega =\Omega (t^{n})$ and the weak quasi-homogeneous
pair $(h,\tilde{h})$ is also regular,
then the associated weak
$\mathcal F$-manifold is $(M_{0},\bullet ,\tilde{h}, E)$ and is an
$\mathcal F$-manifold.

\item Let $\Omega (t)=(ct^{n}+d)^{-1}$, for
$c,d$ constants.
The pair $(h,\tilde{h})$ together with $E$ is weak quasi-homogeneous on
$\mathbb{H}^{n}\otimes\mathbb{C}/W$ (when $(cd)<0$)
and on $\mathbb{C}^{n}/W$ (when $(cd)>0$). It is regular on the open subset
where
$$
t^{n}\neq \frac{d}{c},\quad t^{n}\neq\frac{(1-d_{1})d}{(1+d_{1})c}.
$$
Moreover, $\tilde{h}$ is flat and $h$ has constant sectional
curvature $4(cd).$

\end{enumerate}

\end{thm}

\begin{proof}
The first statement follows from Proposition \ref{comp}.
The second statement is
a consequence of Proposition \ref{f}:
note that $g(E)=\tilde{g}\left(\frac{\partial}{\partial t^{1}}\right)
=dt^{n}$. The third statement uses the proof of
Proposition \ref{abra}. Note that
the endomorphism $T$ from Definition \ref{defini}
has the following expression:
$$
T(u)=\sum_{i=1}^{n}\left( (d_{i}-1)u_{i}+\frac{cu_{1}}{ct^{n}+d}d_{n-i+1}
t^{n-i+1}\right) dt^{i}-\frac{cu(E)}{ct^{n}+d}dt^{n},
$$
for every $u=\sum_{i=1}^{n}u_{i}dt^{i}.$
The regularity condition can be easily checked.
\end{proof}

\medskip

{}From the symmetries of the tensors and the flatness of the metric,
one may integrate the equations, via the Poincar\'e lemma, and
express the tensors as derivatives, with respect to the flat
coordinates, of a scalar prepotential. The differential equation
satisfied by the prepotential being the celebrated
Witten-Dijkgraff-Verlinde-Verlinde (or WDVV) equation.

\medskip

Thus given a Frobenius manifold with prepotential $F$ one may
conformally rescale the metrics, derive new flat coordinates and
multiplication, and calculate the new prepotential ${\tilde F}\,.$
\begin{eqnarray*}
F & \rightarrow & \{ \bullet_g\,,{\tilde g} \} \\
\downarrow & & \quad\downarrow \\
{\tilde F} & \leftarrow &\{ \bullet_h\,,{\tilde h} \}
\end{eqnarray*}
This gives rise to an $SL(2,\mathbb{C})$-symmetry on solution space of the WDVV
equation.

\medskip

\noindent{\bf Example:}
\medskip

\noindent Starting with the prepotential\footnote{For notational convenience indices
are dropped in these example {\sl only}, so $t_i=t^i\,.$}

\[
F = \frac{1}{2} t_1^2 t_3 + \frac{1}{2} t_1 t_2^2 + f(t_2,t_3)
\]
where $f$ satisfies the differential equation
\[
f_{333}=f_{223}^2 - f_{233} f_{222}
\]
one obtains the new solution

\[
{\tilde F} = \frac{1}{2} {\t}_1^2 t_3 + \frac{1}{2} \t_1 \t_2^2 +\left\{\frac{c \t_2^4}{8(c\t_3+d)} + (c\t_3+d)^2 f\left(
\frac{\t_2}{c \t_3+d},\frac{a\t_3+b}{c\t_3+d}\right)\right\}
\]
where $ad-bc=1\,.$ Note that this is a transformation on solutions on the WDVV equation: the transformation
breaks the linearity condition on the Euler vector field (except in the very special case identified in \cite{D})
and so does not generate new examples of Frobenius manifolds.

\medskip

As mentioned in the introduction, these conformally flat pencils will automatically generate
bi-Hamiltonian structures and hence certain integrable hierarchies of evolution equations.
The properties of these hierarchies will be considered elsewhere.

\vskip 1cm

\noindent{Authors' addresses:}\par

\bigskip

\begin{tabular}{ll}
Liana David: & Institute of Mathematics of the
Romanian Academy, \\&Calea Grivitei nr 21, Bucharest, Romania;\\&
e-mail: liana.david@imar.ro; lili@mail.dnttm.ro\\&\\

&\\
Ian A.B. Strachan: &Department of Mathematics, University of
Glasgow, \\&Glasgow G12 8QW, U.K.; \\&e-mail: i.strachan@maths.gla.ac.uk\\

\end{tabular}

\end{document}